\documentclass[a4paper,twoside,10pt,fleqn]{article}

\usepackage{graphics}
\usepackage{latexsym}
\usepackage{amsmath}
\usepackage{amssymb}
\usepackage{epsf}

\newcommand{\abc}{
\renewcommand{\theenumi}{\alph{enumi}}
\renewcommand{\labelenumi}{\theenumi)}
\itemsep 0pt}

\newcommand{\commentaar}[1]{{}}

\newcommand{\inprod}[2]{\langle #1, #2 \rangle}
\newcommand{\partieel}[2]{\frac{\partial #1}{\partial #2}}

\newcommand{\Diff}{{\rm Diff}}
\newcommand{\Hess}{\textrm{Hess}}
\newcommand{\calO}{{\cal O}}
\newcommand{\cinf}{C^{\infty}}
\newcommand{\eps}{\varepsilon}

\newcommand{\fN}{\mathbb{N}}
\newcommand{\fR}{\mathbb{R}}
\newcommand{\grad}{\textrm{grad}}

\newcommand{\rand}{\partial}
\newcommand{\vol}{{\rm vol}}

\renewcommand{\empty}{\varnothing}

\newtheorem{theorem}{Theorem}
\newtheorem{definition}{Definition}
\newtheorem{proposition}[theorem]{Proposition}
\newtheorem{lemma}[theorem]{Lemma}

\newcommand{\versie}{05/dec/2007}
\newcommand{\kortetitel}{Differentiability}
\newcommand{\korteauteur}{Hoveijn}

\title{Differentiability of the volume of a region enclosed by level sets}
\author{I. Hoveijn\\University of Groningen, Department of Mathematics\\
P.O. Box 800, 9700 AV Groningen, The Netherlands}
\date{\versie}

\begin{document}
\maketitle

\begin{abstract}\noindent
The level of a function $f$ on $\fR^n$ encloses a region. The volume
of a region between two such levels depends on both levels. Fixing one
of them the volume becomes a function of the remaining level. We show
that if the function $f$ is smooth, the volume function is again
smooth for regular values of $f$. For critical values of $f$ the
volume function is only finitely differentiable. The initial
motivation for this study comes from Radiotherapy, where such volume
functions are used in an optimization process. Thus their
differentiability properties become important.
\end{abstract}

\section{Introduction}\label{sec:intro}

The volume of a set enclosed by two different level sets of a function
$f: \fR^n \to \fR$ depends on both levels. Here we fix one level and
the question we want to address is the differentiability of the volume
as a function of the varying level for arbitrary dimension $n$. It
will turn out that under mild conditions, among which smoothness of
$f$, this function is again smooth for all regular values of $f$, but
only finitely differentiable at critical values of $f$. Moreover we
also consider $f$ on a compact subset $V$ of its domain and
subsequently pose the same question for level sets restricted to
$V$. Then we get a similar differentiability result when we also
include the levels of $f$ restricted to the boundary of $V$. In
section \ref{sec:stats} we will more precisely define the situation we
consider and give precise statements of our results.

The initial motivation for this study comes from Radiotherapy. A
patient is treated with ionizing radiation causing energy release per
unit mass or volume inside the patient, which is called dose. The
above mentioned function $f$ represents this dose and the set $V$
represents a patient's organ or a tumour region. The therapeutic
outcome of Radiotherapy treatment not only depends on the dose on the
tumour but also on the dose on healthy organs. Since it is usually
hard to estimate the effect of dose in three dimensons, one reduces it
for each tumour and organ to a so called \emph{dose-volume histogram},
see \cite{khfm}. In it one records for each dose value which fraction
of the volume receives at least this dose. The resulting collection of
dose-volume pairs is called the \emph{cumulative} dose-volume
histogram, which is closely related to the volume function introduced
in the first paragraph. It is well established that the larger the
fraction of the volume of a tumour receiving a prescribed dose, the
larger the probability that the tumour is eradicated. For healthy
organs and tissues the situation is less clear but some evidence
exists that damage can be estimated from dose-volume histograms, see
\cite{khfm} and references therein.

In Radiotherapy treatment planning one tries to find an optimum for a
sufficiently high dose on a tumour and a sufficiently low dose on
healthy organs. The object function of this optimization process
depends in particular on certain dose-volume pairs $(d,v(d))$ on the
graph of the cumulative dose-volume histogram, see \cite{wm}. In the
optimization process we have a family of dose functions and thus a
family of dose-volume histograms, parametrized by the optimization
variables. Anticipating the result on finite differentiability of the
dose-volume histogram at a critical value of the dose function, we
conclude that the object function is finitely differentiable at values
of the optimization variables for which $d$ is a critical value of the
dose function. Because many numerical methods to solve such problems
assume differentiability of the object function to some order,
differentiability of the dose-volume function becomes important.

The remaining part of this article is organized as follows. In section
\ref{sec:stats} we give a summary of the results, preceded by
definitions. In section \ref{sec:proofs} we sketch the idea of the
proof and the results are proved in a series of propositions. Some of
the more elaborate computations are summarized in an appendix.

\section{Statement of results}\label{sec:stats}
%

\subsection{Definitions}\label{sec:defs}
Our main object is a smooth positive function $f$ which is bounded and
whose levelsets are compact. These properties are sufficient for our
results, we do not claim necessity. We now define a function class for
future reference.

\begin{definition}\label{def:functieklasse}
Let $C$ be the class of functions $f: \fR^n \to \fR$ satisfying the following\\[-5ex]
\begin{enumerate}\abc
\item $f$ is $\cinf$,
\item $f$ is positive,
\item $f$ is decreasing, that is for each $\eps > 0$ there is a compact
$K \subset \fR^n$ such that for all $x \in \fR^n \setminus K$, $f(x) < \eps$.
\end{enumerate}
\end{definition}

We define the levelsets $N_c$ of $f \in C$ and sets $U_c$ enclosed by
levels in a straightforward manner.
\begin{eqnarray*}
N_c &=& \{x \in \fR^n \;|\; f(x) = c \},\\
U_c &=& \{x \in \fR^n \;|\; f(x) \geq c \}.
\end{eqnarray*}
We wish to study the study the differentiability of the volume
$\vol_n(U_c)$ of $U_c$ as a function of $c$. Here $\vol_n(A)$ is the
standard volume (Lebesgue measure) in $\fR^n$ of a measurable set $A
\subset \fR^n$
\begin{displaymath}
\vol_n(A) = \int_A 1\;dx.
\end{displaymath}

When we restrict $f$ to a subset $V$ of $\fR^n$ we make similar
definitions. To be more precise, let $V \subset \fR^n$ be a compact,
connected, $n$-dimensional $\cinf$-manifold with boundary $\rand V$
which is a $(n-1)$-dimensional $\cinf$-manifold. The restriction of
$f$ to $V$ will be denoted by $f|_V$. We define $V_c = U_c \cap V$ and
the volume function we consider for this case is $\vol(V_c)$.

As already indicated in the Introduction, the function $\vol$ need not
be continuous. Indeed, let $f$ be constant on an open neighbourhood of
$V$ then it is easily seen that $\vol$ is similar, upon scaling and
translation, to the Heaviside function. In order to avoid this and
other degeneracies we assume the following non-degeneracy conditions
on $f$, $f|_V$ and $V$.

\begin{definition}\label{def:nondeg} Non-degeneracy conditions.\\[-5ex]
\begin{enumerate}\abc \parskip 0pt \itemsep 0pt
\item \label{itm:isolcrit} Critical points of $f$ are non-degenerate, that is if $x$
is a critical point of $f$ then $\det(\Hess\;u(x)) \neq 0$.
\item \label{itm:inwcrit} If $x$ is a critical point of $f$ then $x \not \in \rand V$.
\item \label{itm:isolcrit2} Critical points of $f|_V$ are non-degenerate.
\item \label{itm:fine} $f$ is \emph{fine}, that is if $x$ and
$y$ are critical points and $x \neq y$ then $f(x) \neq f(y)$.
\item \label{itm:fine2} $f|_V$ is fine.
\end{enumerate}
\end{definition}

\textbf{Remarks}
\begin{enumerate}
\item Conditions \ref{itm:isolcrit}, \ref{itm:inwcrit} and
\ref{itm:isolcrit2} are essential for our proofs. If one of them is not
satisfied our standard forms, see section \ref{sec:proofs}, of $f$
and $V$ are no longer valid. In order to proceed we would need
higher order information about $f$ and $V$. Moreover the critical
point would not be stable under small perturbations. At the moment we
are not interested in such a situation. Conditions \ref{itm:fine} and
\ref{itm:fine2} are for convenience only and can easily be dropped.
\item In principle, the sets $U_h$, $V$ and $V_h$ need not be
connected. If one of them consists of several components, the
construction in the following sections can be carried out for each
component separately. Therefore without loss of generality we may as
well assume connectedness.
\item In Radiotherapy one sometimes uses the so called
\emph{differential} dose-volume histogram, see \cite{khfm}. However in
general this is not a function. In terms of our function $f$ it is in
fact the image measure of the standard measure on $\fR^3$ under $f$ on
$\fR$, see \cite{cohn}. Instead of studying the volume function via
$f$ as defined in this section one could also study the volume
function via this measure. The present approach however seems to be
simpler.
\end{enumerate}

\subsection{Results}\label{sec:results}
In order to state the results we will make a distinction between
regular values and critical values of $f|_V$. In case of a critical
value of $f|_V$ we make a further distinction whether the critical
point is in $\rand V$ or not. Results for $f$ and $f|_V$ are identical
when the critical point is not in $\rand V$, therefore they are not
stated separately. Our proofs are valid only for $f$ in class $C$, see
definition \ref{def:functieklasse}, satisfying the nondegeneracy conditions
in definition \ref{def:nondeg}.

The first theorem states that for all regular values of the function
$f$, which means for almost all values by virtue of the nondegeneracy
conditions, the volume function is a smooth function of the level.

\begin{theorem}\label{the:reg}
Let $0 \in \fR$ be a regular value of $f|_V$ then $\vol(V_h)$ is a
smooth function of $h$ at $0$.
\end{theorem}

The second and third theorem state that at a critical value the
differentiability of the volume function is finite. The order of
differentiability depends on the dimension of the domain of $f$. The
nature of the discontinuity depends on the Morse index of $f$ at the
critical point. Together with the proofs we give details about the
discontinuity in the next section.

\begin{theorem}\label{the:crit1}
Let $0 \in \fR$ be a critical value of $f|_V$ and let $0 \not \in
\rand V$ be the critical point. Then the $\lceil\frac{n}{2}\rceil$-th derivative
of $\vol(V_h)$ is discontinuous at $h=0$, all lower order derivatives are continuous.
\end{theorem}

\begin{theorem}\label{the:crit2}
Let $0 \in \fR$ be a critical value of $f|_{\rand V}$ and let $0 \in
\rand V$ be the critical point. Then the $\lceil\frac{n+1}{2}\rceil$-th derivative
of $\vol(V_h)$ is discontinuous at $h=0$, all lower order derivatives are continuous.
\end{theorem}

In a Radiotherapy setting, which originally motivated this study, the
dimension of the domain of $f$, representing dose, is 3. This means
that at a critical point of dose, the volume function is not even
twice continuously differentiable. As mentioned in the introduction
one aimes at finding an optimum for a suffiently high dose on the
tumour and a sufficiently low dose on healthy organs. The objective
function in this optimization process depends on the volume
function. Many iterative optimization methods use a quasi-Newton
method in the background and thus require differentiability to second
order. Such a method is not guaranteed to be well behaved near a
critical value of dose.

\section{Proof of results}\label{sec:proofs}

We will work in the class of $\cinf$-functions. This class is closed
under the action of the group of $\cinf$-transformations. Therefore we
have the notion of $\cinf$-equivalence of functions. We use this to
put the function at hand into a suitable standard form. In general
however, this standard form is only valid locally in a small open ball.

Now our aim is to compute the volume of the set enclosed by two
levelsets of a function, which is a global rather than a local
problem. A natural way to look at this set is by ``sweeping out''
using the gradient flow of the function. Using compactness we can then
turn this into a local problem considering flow boxes of the gradient
flow starting in small subsets, here we take $(n-1)$-simplices, of the
levelset corresponding to the lower level and ending in the levelset
corresponding to the higher level. By a suitable
$\cinf$-transformation we turn each flow box into a Cartesian product
of a $(n-1)$-simplex and an interval. This greatly simplifies finding
the volume of the enclosed set and its dependence on the level.

A complication in this procedure is that a general
$\cinf$-transformation does not map a pair $(f,\grad f)$ into a new
pair $(g,\grad g)$. In order to achieve the latter, the transformation
would have to preserve the inner product which is used to define
$\grad$. For our purposes this is too stringent a restriction. Another
reason not to use the gradient flow is the following. We also wish to
consider the set enclosed by levelsets restricted to a set $V$. The
gradient flow is not necessarily tangent to $\rand V$, thus the flow
box we construct might not be restricted to $V$. Instead we use the
flow of a $\cinf$-vectorfield which is only transversal to the
levelsets and tangent to $\rand V$ if necessary. These properties are
preserved by a general $\cinf$-transformation.

The tansformations we apply to a flow box do not in general preserve
its volume. Here, however, we are only interested in the $h$
dependence of this volume, not in its numerical value. Therefore we
may apply affine transformations without any further
considerations. With other transformations we have to be more careful
and we will take them into account at the appropriate places.

Let us sketch the steps in the proofs of theorems \ref{the:reg},
\ref{the:crit1} and \ref{the:crit2}. If $0$ is a regular value of $f$,
all points in $V_h$ are regular for $h$ small enough. We first
construct a finite number of boxes $B_i$ covering $V_h$ using a
triangulation of $N_0$ and a regular flow from $N_0$ to
$N_h$. Assuming the boxes $B_i$ are small enough we put each of them
in standard form by several local $\cinf$-transformations. The first
transformation parallellizes the flow from $N_0$ to $N_h$. The second
transformation is linear and preserves the parallellity of the flow
but makes it perpendicular to $N_0$ at $0$ and parallel to the last
basis vector $e_n$ of $\fR^n$. The third transform takes $f$ into a
local standard form preserving all of the previous. The result is that
$B_i$ is transformed to the Cartesian product of a simplex in $N_0$
and the interval $[0,h]$. The conclusion is that the volume of $B_i$
is a smooth function of $h$.

If $0$ is a critical value of $f$ we use the fact that $0$ is the only
critical point on $N_0$. We construct one special box $B_0$ containing
$0$ and away from $0$ we use the same construction as
above. Differentiability is then determined by $\vol(B_0)$. Again we
use a transformation that takes $f$ into standard form, but now at the
critical point.

The main part of the proof of theorem \ref{the:reg} is the
construction of the boxes and putting them into a standard form. In
the proofs of the other theorems the emphasis is on computing the
volume of box $B_0$.

Where necessary we assume the existence of a standard basis and a
standard inner product.

\subsection{Proof of theorem \ref{the:reg}}\label{sec:thereg}

\begin{proposition}\label{pro:doosjes}
Let $0$ be a regular value of $f$, then for sufficiently small $h>0$
there exists a finite collection of sets $B_i$ with $i \in I \subset
\fN$ satisfying
\begin{enumerate}\abc \parskip 0pt \itemsep 0pt
\item $V_h = \cup_{i \in I} B_i$
\item $\vol_n(V_h) = \sum_{i \in I} \vol_n(B_i)$.
\end{enumerate}
\end{proposition}

\textbf{Proof of proposition \ref{pro:doosjes}.}\\
\textbf{Construction.} If $0$ is a regular value of $f$ then $N_0$ is a
smooth manifold. The non-degeneracy conditions imply that critical
values of $f$ are isolated, therefore an $h_0>0$ exists such that all
$h \in [0,h_0]$ are regular values. Then all $N_h$ are diffeomorphic
to $N_0$, see \cite{mil}. \commentaar{(met het voorgaande doen we
niets?)} First we assume that $N_0$ does not intersect the boundary
$\rand V$ of $V$. Since $N_0$ is a compact $\cinf$-manifold it allows
a finite triangulation with $(n-1)$-simplices $\sigma_{i,n-1}$ and $i
\in I \subset \fN$, see \cite{spi}. One of the properties of a
triangulation is that for $i \neq j$ either $\sigma_{i,n-1} \cap
\sigma_{j,n-1} = \empty$ or $\sigma_{i,n-1} \cap \sigma_{j,n-1} =
\sigma_{k,n-2}$ for some $k$ and $(n-2)$-simplex $\sigma_{k,n-2}$. Let
$X$ be a $\cinf$-vector field transversal to $N_h$ for all $h \in
[0,h_0]$. More precisely we impose the condition that there is an
$\eps > 0$ such that $|\inprod{n_h(x)}{X(x)}| > \eps$ for all $x \in
N_h$ and $h \in [0,h_0]$, where $n_h(x)$ is a unit normal to $N_h$ at
$x$. Such a vectorfield exists, for example $\grad f$. If $N_0$
intersects $\rand V$ then the non-degeneracy conditions imply that the
intersection of $N_h$ and $\rand V$ is transverse for all $h \in
[0,h_0]$. We restrict to $N_0 \cap V$ which is still a compact
$\cinf$-manifold. Now we impose one more condition on the vector field
$X$ namely that it is tangent to $\rand V \cap V_{h_0}$. Let $\Phi \in
\Diff(\fR^n)$ be the flow of $X$ with $\Phi(x,0)=x$. Finally we define
\begin{displaymath}
B_i = \{\Phi(x,t) \;|\; x \in \sigma_{i,n-1}, t \in [0,T] \} \cap V_h.
\end{displaymath}
\textbf{Proof of {\boldmath $V_h = \cup_{i \in I} B_i$}.} Clearly $\cup_{i \in I}
B_i \subset V_h$. Suppose $x \in V_h$, then $\frac{d}{dt}
f(\Phi(x,t))$ is strictly increasing or decreasing since
$|\inprod{n_h(x)}{X(x)}| > \eps$. In either case a finite $t_0$ exists
such that $f(\Phi(x,t_0)) = 0$. This means $x_0 =
\Phi(x,t_0) \in N_0$, therefore an $i$ exists such that $x_0 \in
\sigma_{i,n-1}$ which implies $x \in B_i$. Compactness of $V_h$
guarantees that a $T>0$ exists such that for all $x \in V_h$, $|t_0|
\in [0,T]$. The conclusion is that $V_h = \cup_{i \in I} B_i$.\\
\textbf{Proof of {\boldmath $\vol_n(V_h) = \sum_{i \in I} \vol_n(B_i)$}.} It
suffices to show that for $i \neq j$, $\vol_n(B_i \cap B_j) =
0$. Suppose $B_i \cap B_j \neq \empty$ and $x \in B_i \cap B_j$, then
there is a $t_0$ such that $x_0 = \Phi(x,t_0) \in N_0$. By definition
$\Phi(x,t) \in B_i \cap B_j$ therefore $x_0 \in B_i \cap B_j$ which
means $x_0 \in \sigma_{i,n-1} \cap \sigma_{j,n-1} =
\sigma_{k,n-2}$. From this we conclude $B_i \cap B_j = \{\Phi(x,t)
\;|\; x \in \sigma_{k,n-2}, t \in [0,T] \} \cap V_h$. But then
$\dim(B_i \cap B_j) = n-1$ and therefore $\vol_n(B_i \cap B_j) = 0$.
\hfill $\blacksquare$

The next step is to put the boxes $B_i$ of proposition
\ref{pro:doosjes} into a standard form. To do this we also need a
local standard form of the function $f$.

\begin{proposition}\label{pro:fnfs}
Let $f$ be a function as in definition \ref{def:functieklasse}, satisfying
the non-degeneracy conditions in definition \ref{def:nondeg} and
$f(0)=0$. We distinguish three different cases:
\begin{enumerate}\abc \parskip 0pt \itemsep 5pt
\item $0$ is a regular point of $f$,
\item $0$ is a critical point of $f$,
\item $0$ is a critical point of $f|_{\rand V}$.
\end{enumerate}
Then an open neighbourhood $\calO$ of $0$ and a diffeomorphism $\Phi$
exist such that $F = \Phi_*f$ takes one of the forms:
\begin{enumerate}\abc \parskip 0pt \itemsep 5pt
\item $F(\xi, \eta) = \eta$, with $(\xi, \eta) \in (\fR^{n-1} \times
\fR) \cap \calO$,
\item $F(\xi, \eta) = \sum_{i=1}^p \xi_i^2 - \sum_{i=1}^p \eta_i^2$, with $(\xi,
\eta) \in (\fR^p \times \fR^q) \cap \calO$ and $p+q=n$,
\item $F(\xi, \eta, \zeta) = \zeta + \sum_{i=1}^p \xi_i^2 - \sum_{i=1}^p
\eta_i^2$ and $\rand V$ is given by $\zeta=0$ with $(\xi, \eta, \zeta)
\in (\fR^p \times \fR^q \times \fR) \cap \calO$ and $p+q=n-1$.
\end{enumerate}
\end{proposition}

\textbf{Remarks}
\begin{enumerate}
\item Case b) of proposition \ref{pro:fnfs} is called the Morse lemma. A
critical point in this case has \textit{Morse index} $q$, but some
times it is more convenient say it has \textit{Morse type} $(p,q)$.
\item Due to non-degenracy condition 2, $0$ is a critical point of
$f|_{\rand V}$ as soon as $0$ is a critical point of $f|_V$.
\end{enumerate}

In the proof of \ref{pro:fnfs} we will need a lemma which we only
state here, for a proof see \cite{mil}.

\begin{lemma}\label{lem:hetlemma}
Let $f$ be a $\cinf$ function on a convex neighbourhood $\calO \subset
\fR^n$ of $0$ with $f(0)=0$. Then $f(x) = \sum_{i=1}^n x_i f_i(x)$ for
certain $\cinf$ functions $f_i$ with $f_i(0) = \partieel{}{x_i}f(0)$.
\end{lemma}

\textbf{Proof of proposition \ref{pro:fnfs}.}
\begin{enumerate}\abc \parskip 0pt \itemsep 5pt
\item Since $0$ is a regular point of $f$ there is a nonzero vector $a$
such that $\grad f(0) = a$. After an orthogonal transformation we may
assume that with respect to coordinates $(x,y) \in \fR^{n-1} \times
\fR$ we have $\grad f(0,0) = (0,|a|)$. Now we define new coordinates
by the diffeomorphism $\Phi: (x,y) \mapsto (x,f(x,y))$, then $F =
\Phi_*f$ takes the desired form.
\item See \cite{mil}.
\item If $0$ is a critical point of $f_V$, the the tangent spaces of
$N_0$ and $V$ at $0$ coincide. By the non-degeneracy conditions $0$ is
a regular point of $f$. Using the arguments of case a) we assume that
we already transformed to coordinates $(x,y) \in \fR^{n-1} \times \fR$
such that $\grad f(0,0) = (0,|a|)$. Furthermore we made the assumption
that $\rand V$ is a smooth manifold, so at least locally it is the level
set of a smooth function $g$. Now we apply part a) to bring $g$ into
standard form, then on new coordinates $(u,v) \in \fR^{n-1} \times
\fR$, $\rand V$ is locally given by $v=0$. After scaling in the $v$
direction the function $f$ satisfies: $f(0,0) = 0$,
$\partieel{}{u_i}f(0,0) = 0$ for $i \in \{1,\ldots,n-1\}$ and
$\partieel{}{v}f(0,0)=1$. The remainder of the proof is only a slight
adaption of the proof in \cite{mil} for case b), but included here for the
sake of completeness. Applying lemma \ref{lem:hetlemma} to $f$ and its
partial derivatives we get
\begin{displaymath}
f(u,v) = \sum_{i,j < n} u_i u_j \alpha_{ij}(u,v) + v f_n(u,v)
\end{displaymath}
where $\alpha_{ij}$ and $f_n$ are smooth functions and
$f_n(0,0)=1$. As a first step we apply the transformation $(u,v)
\mapsto (u, v f_n(u,v)) = (u,z)$, preserving the standard form of $g$
since $f_n(0,0) \neq 0$. Then $f$ takes the form
\begin{displaymath}
f(u,z) = \sum_{i,j < n} u_i u_j \alpha_{ij}(u,z) + z
\end{displaymath}
with new functions $\alpha_{ij}$. From now on we only apply
transformations of the form $(u,z) \mapsto (\phi(u,z),z)$ and we
proceed by induction. We assume that
\begin{displaymath}
f(u,z) = \pm u^2_1 \cdots \pm u^2_{k-1} + \sum_{i,j = k}^{n-1} u_i u_j
\alpha_{ij}(u,z) + z
\end{displaymath}
for a certain $k>0$. Now let $\hat{u}_i = u_i$ for $i \neq k$ and
$\hat{z} = z$ and
\begin{displaymath}
\hat{u}_k = \sqrt{\alpha_{kk}(u,z)}(u_k + \sum_{i=k+1}^{n-1} u_i
\alpha_{ik}(u,z)/\alpha_{kk}(u,z))
\end{displaymath}
In order to define $\hat{u}_k$ it may be necessary to permute the rows
of $\alpha_{ij}$ so that $\alpha_{kk}(u,z) \neq 0$. Such a permutation
exists because $\det(\Hess f(0,0)) \neq 0$ which means there is at
least one $i \in \{k,\ldots,n-1\}$ such that $\alpha_{ik}(0,0) \neq
0$. Then by continuity there is a neighbourhood of $(0,0)$ such that
$\alpha_{ik}(u,z) \neq 0$. This means that in each induction step the
neighbourhood on which our result holds might shrink. Since we only
need a finitie number of steps this does not cause any
problems. Dropping the hats we get in new coordinates
\begin{displaymath}
f(u,z) = \pm u^2_1 \cdots \pm u^2_k + \sum_{i,j = k+1}^{n-1} u_i u_j \alpha_{ij}(u,z) + z
\end{displaymath}
Renaming the variables we arrive at the desired form of $F = \Phi_*f$,
where $\Phi$ is the composition of the transformations in each
induction step.\hfill $\blacksquare$
\end{enumerate}

\begin{proposition}\label{pro:doosjesnf}
Let $\{B_i\}_{i \in I}$ be a collection of boxes as in proposition
\ref{pro:doosjes}. Let $x$ be a point in the interior of
$\sigma_{i,n-1}$ for some $i \in I$, without loss of generality we
assume that $x=0$. Then a diffeomorphism $\Phi$ exists such that
\begin{displaymath}
\vol_n(B_i) = \int_{B_i} 1\,dx = \int_0^h \left[ \int_{\sigma_{i,n-1}}
J_{\Phi} \,d\xi \right] \,d\eta,
\end{displaymath}
where $J_{\Phi}$ is the Jacobian of $\Phi$. Moreover $\vol_n(B_i)$ is
a smooth function of $h$.
\end{proposition}

\textbf{Proof of proposition \ref{pro:doosjesnf}.}
The vector field in the construction of box $B_i$ has no stationary
point therefore it can be parallellized by a diffeomorphism $\Phi_1$,
see \cite{spi}. The next two transformations preserve parallellity
because they are linear. By a linear diffeomorphism $\Phi_2$ we can
arrange that $X$ is perpendicular to $N_0$ in $0$. By another linear
diffeomorphism $\Phi_3$ we rotate such that $X$ is parallel to the
last basis vector of $\fR^n$. The last diffeomorphism $\Phi_4 : (x,y)
\mapsto (x,f(x,y))$ takes the function $f$ into local standard form,
see proposition \ref{pro:fnfs}. Since $\Phi_4$ is a position dependent
shift in the direction of the vectorfield, parallellity is again
preserved. However the parametrization of the integral curves will
change in general. Then $\Phi = \Phi_4 \circ \Phi_3
\circ \Phi_2 \circ \Phi_1$ is again a diffeomorphism and on new
coordinates $(\xi,\eta)$ we have $B_i = \sigma_{i,n-1} \times
[0,h]$. Since $J_{\Phi}$ is a smooth function and $\vol_n(B_i)$
depends on $h$ only via the upper limit of the outer integral,
$\vol_n(B_i)$ is a smooth function of $h$.\hfill $\blacksquare$

Using the previous propositions we are able to prove theorem \ref{the:reg}.

\textbf{Proof of theorem \ref{the:reg}.}
Construct boxes $B_i$ as in proposition \ref{pro:doosjes}. Then by the
same proposition $\vol_n(V_h) = \sum_{i \in I} \vol_n(B_i)$. In the
latter $\vol_n(B_i)$ is a smooth function of $h$ by proposition
\ref{pro:doosjesnf}. The finite number of boxes guarantees that also
$\vol_n(V_h)$ is a smooth function of $h$. \hfill $\blacksquare$

\subsection{Proof of theorem \ref{the:crit1}}\label{sec:thecrit1}
In the previous section all points were regular points of $f$. Here
too points will be regular except the point $0$ which now is a
critical point. Therefore we will use the same construction of boxes
as in the previous section. Only the box $B_0$ containing $0$ will be
treated differently. This means that differentiability in this
situation is determined by $\vol_n(B_0)$.

A critical point with Morse index 0 is a minimum of $f$. If the Morse
index is $n$ the critical point is a maximum of $f$. Differentiability
for minima and maxima is very similar so we only state a result for
one of them. A critical point with Morse index $q$ where $0 < q < n$ is
called a saddle.

\begin{proposition}\label{pro:min}
Let $0$ be a non-degenerate critical point of $f$ with Morse index
$0$, so $0$ is a local minimum. Then the $\lceil\frac{n}{2}\rceil$-th
derivative of $\vol(V_h)$ is not continuous at $h=0$, all lower order
derivatives are continuous.
\end{proposition}

\textbf{Proof of proposition \ref{pro:min}.}
Let $h \geq 0$, then $0$ is the minimal value of $f$ so $V_{-h}$ is
empty. The level set $N_0$ only contains $0$. Using the neighbourhood
$\calO$ of proposition \ref{pro:fnfs} on which $f$ takes its standard
form we compute $\vol_n(V_h)$ using polar coordinates $(r, \varphi)$,
then that $u(r, \varphi) = r^2$. Here we assume that $h$ is small
enough so that $V_h \subset \calO$. The Jacobian of the transformation
$\Phi$ in proposition \ref{pro:fnfs} will be denoted by $J_{\Phi}$ and
the Jacobian of changing to polar coordinates by $r^{n-1}
g_n(\varphi)$. Since $\Phi$ is non-singular and $\cinf$ we can split
$J_{\Phi} = c + K_{\Phi}$ where $c \neq 0$ and both $J_{\Phi}$ and
$K_{\Phi}$ are $\cinf$. Then we have
\begin{displaymath}
\vol_n(V_h) = \int_{V_h} 1\,dx = \int_{V_h} J_{\Phi} \,d\xi = c \int
\int_0^{\sqrt{h}} r^{n-1} \,dr\,g_n(\varphi)\,d\varphi + \int_{V_h}
K_{\Phi} \,d\xi.
\end{displaymath}
For a non-differentiability result it is enough to consider the first
integral in the last expression, since this integral contains the
lowest order terms in $h$ and we are interested in $h \to 0$ only. Let
$a_n$ be the ``area'' of $S^{n-1}$ then
\begin{displaymath}
c \int \int_0^{\sqrt{h}} r^{n-1} \,dr\,g_n(\varphi)\,d\varphi = c\,a_n
\int_0^{\sqrt{h}} r^{n-1} \,dr = \frac{c a_n}{n} h^{n/2}.
\end{displaymath}
Thus we obtain the result that the $\lceil\frac{n}{2}\rceil$-th derivative is
discontinuous at $h=0$.\hfill $\blacksquare$

\begin{proposition}\label{pro:zad}
Let $0$ be a non-degenerate critical point of $f$. Assume the Morse
type of $0$ is $(p,q)$ with $p \neq 0$ and $q \neq 0$. Then the
$\lceil\frac{n}{2}\rceil$-th derivative of $\vol(V_h)$ is
discontinuous at $h=0$. For both $p$ and $q$ even the discontinuity is
a jump, for $p$ and $q$ odd it is a log-like singularity and for $p+q$
odd it is a root-like singularity. All lower order derivatives are
continuous.
\end{proposition}

\textbf{Proof of proposition \ref{pro:zad}.}
Let $h_0$ be small enough so that $0$ is the only critical point of $f$ in
$V_{h_0}$. Furthermore let $\calO$ be a neighbourhood as in proposition
\ref{pro:fnfs} such that $f$ can be put into standard form
c). Changing to cylinder coordinates $(r, \varphi, s, \psi)$ we may
assume that $f(r, \varphi, s, \psi) = r^2 - s^2$. Then an $\eps > 0$
exists such that $B = \{(r, \varphi, s, \psi) \in \fR \times S^{p-1}
\times \fR \times S^{q-1} \;|\; r+s \leq \eps\}$ is a subset of $\calO$. The
intersection of $\rand B$ and $N_0$ is transversal so again taking
$h_0$ small enough we may assume that $\rand B$ transversally
intersects $N_h$ for all $h \in [0,h_0]$. Then $V_h \setminus B$ is a
compact $\cinf$-manifold (with boundary) so we can apply the
construction of boxes as in proposition \ref{pro:doosjes}. Where we
have to impose the additional condition that the vector field $X$ is
tangent to the boundary $\rand B$ of $B$. By proposition
\ref{pro:doosjesnf} $\sum_{i \in I} \vol_n(B_i)$ is a smooth function of
$h$. If we now set $B_0 = B \cap V_h$ then $\vol_n(V_h) =
\sum_{i \in I} \vol_n(B_i) + \vol_n(B_0)$. Thus differentiability of
$\vol_n(V_h)$ is determined by the differentiability of $\vol_n(B_0)$.

The Jacobian of the transformation $\Phi$ in proposition
\ref{pro:fnfs} will be denoted by $J_{\Phi}$ and the Jacobian of
changing to polar coordinates by $r^{p-1} s^{q-1} g_p(\varphi)
g_q(\psi)$. The last transformation we apply is a scaling so that
$B$ is bounded by $r=0$, $s=0$ and $r+s=1$. Then we get
\begin{displaymath}
\vol_n(B_0) = \int_{B_0} 1\,dx = \int_{B_0} J_{\Phi} \,d\xi\,d\eta =
c \int_{B_0} 1 \,d\xi\,d\eta + \int_{B_0} K_{\Phi} \,d\xi\,d\eta.
\end{displaymath}
For our result we only need to compute the first integral in the last
expression, using the same arguments as in the proof of proposition
\ref{pro:min}. The actual computation can be found in appendix \ref{app:bozad}.
From lemma \ref{lem:nondif} in the same appendix the result on
differentiability follows.\hfill $\blacksquare$

\textbf{Proof of theorem \ref{the:crit1}.} The proof follows from propositions
\ref{pro:min} and \ref{pro:zad}. \hfill $\blacksquare$

\subsection{Proof of theorem \ref{the:crit2}}\label{sec:thecrit2}
In this section $0$ is a critical point of $f|_{\rand V}$, but all
other points are regular. First we consider a critical point which is
a minimum of $f|_{\rand V}$. Since for a maximum we get the same
result we do not state it seperately.

\begin{proposition}\label{pro:relmin}
Let $0$ be a critical point of $f|_{\rand V}$. Assume the Morse type
is $(n-1,0)$ with $n>1$. Then the $\lceil\frac{n+1}{2}\rceil$-st
derivative of $\vol(V_h)$ is not continuous at $h=0$, all lower order
derivatives are continuous.
\end{proposition}

\textbf{Proof of proposition \ref{pro:relmin}.}
The proof is very similar to that of proposition \ref{pro:min}
therefore we will only indicate the essential differences. On a
neighbourhood $\calO$ of $0$ such that $f$ can be put into standard
form c) of proposition \ref{pro:fnfs} we take cylinder coordinates
$(r, \varphi, \zeta) \in \fR \times \fR^{n-2} \times \fR$. Then locally
$f(r, \varphi, \zeta) = r^2 + \zeta$ and $\rand V$ is given by
$\zeta=0$. Splitting the Jacobian as in proposition \ref{pro:min} we obtain
\begin{displaymath}
\vol_n(V_h) = \int_{V_h} 1\,dx = \int_{V_h} J_{\Phi} \,d\xi = c \int
\int_0^{\sqrt{h}} \int_0^{h-r^2} r^{n-2} \,d\zeta\,dr\,g_{n-1}(\varphi)\,d\varphi +
\int_{V_h} K_{\Phi} \,d\xi.
\end{displaymath}
Only evaluating the first integral (cf. proof of proposition
\ref{pro:min}) we get
\begin{displaymath}
c \int \int_0^{\sqrt{h}} \int_0^{h-r^2} r^{n-2} \,d\zeta\,dr\,g_{n-1}(\varphi)\,d\varphi =
2c\frac{a_{n-1}}{n^2-1} h^{(n+1)/2}.
\end{displaymath}
Thus we obtain the result that the $\lceil\frac{n+1}{2}\rceil$-st derivative is
discontinuous at $h=0$.\hfill $\blacksquare$

Next we turn our attention to saddle points of $f|_{\rand V}$.

\begin{proposition}\label{pro:relzad}
Let $0$ be a non-degenerate critical point of $f|_{\rand V}$. Assume
the Morse type of $0$ is $(p,q)$ with $p \neq 0$, $q \neq 0$ and
$p+q=n-1$. Then the $\lceil\frac{n+1}{2}\rceil$-st derivative of
$\vol(V_h)$ is not continuous at $h=0$. For both $p$ and $q$ even the
discontinuity is a jump, for $p$ and $q$ odd it is a log-like
singularity and for $p+q$ odd it is a root-like singularity. All lower
order derivatives are continuous.
\end{proposition}

\textbf{Proof of proposition \ref{pro:relzad}.}
We proceed along the lines of the proof of proposition \ref{pro:zad}
again indicating the main differences only. Let $\calO$ be a
neighbourhood of $0$ such that $f$ can be transformed to standard form
c) of proposition \ref{pro:fnfs}. Taking cylinder coordinates $(r,
\varphi, s, \psi, q\zeta)$ we may assume that $f(r, \varphi, s, \psi,
\zeta) = r^2 - s^2 + \zeta$. Then an $\eps > 0$ exists such that $B =
\{(r, \varphi, s, \psi, \zeta) \in \fR \times S^{p-1} \times \fR \times S^{q-1} \times
\fR \;|\; r+s \leq \eps, \; 0 \leq \zeta \leq \eps \}$ is a subset of
$\calO$. Once again taking $h_0$ small enough we may assume that
$\rand B$ transversally intersects $N_h$ for all $h \in [0,h_0]$. Then
$V_h \setminus B$ is a compact $\cinf$-manifold (with boundary) so
here too we can apply the construction of boxes as in proposition
\ref{pro:doosjes} with the additional condition that the vector field
$X$ is tangent to $\rand B$. We set $B_0 = B \cap V_h$ then
$\vol_n(V_h) = \sum_{i \in I} \vol_n(B_i) + \vol_n(B_0)$. Now
differentiability of $\vol_n(V_h)$ is determined by the
differentiability of $\vol_n(B_0)$. Splitting the Jacobian as in
proposition \ref{pro:zad} we get
\begin{displaymath}
\vol_n(B_0) = \int_{B_0} 1\,dx = \int_{B_0} J_{\Phi} \,d\xi\,d\eta\,d\zeta =
c \int_{B_0} 1 \,d\xi\,d\eta\,d\zeta + \int_{B_0} K_{\Phi} \,d\xi\,d\eta\,d\zeta.
\end{displaymath}
We only compute the first integral in the last expression (cf. proof
of proposition \ref{pro:min}). For the actual computation see appendix
\ref{app:borelzad}. From lemma \ref{lem:nondif2} in the same appendix
the result on differentiability follows.\hfill $\blacksquare$

\textbf{Proof of theorem \ref{the:crit2}.} The proof follows from
propositions \ref{pro:relmin} and \ref{pro:relzad}. \hfill $\blacksquare$

\section*{Acknowledgment}

It is a pleasure to thank Floris Takens for critical reading of an
earlier version of the manuscript and making valuable comments.

\appendix

\section{The volume of $B_0$ containing a saddle point}

\subsection{The volume of $B_0$ in proposition \ref{pro:zad}}\label{app:bozad}

\begin{figure}[htbp]
\setlength{\unitlength}{1mm}
\begin{picture}(100,70)(5,-10)
\put( 20, 5){\epsffile{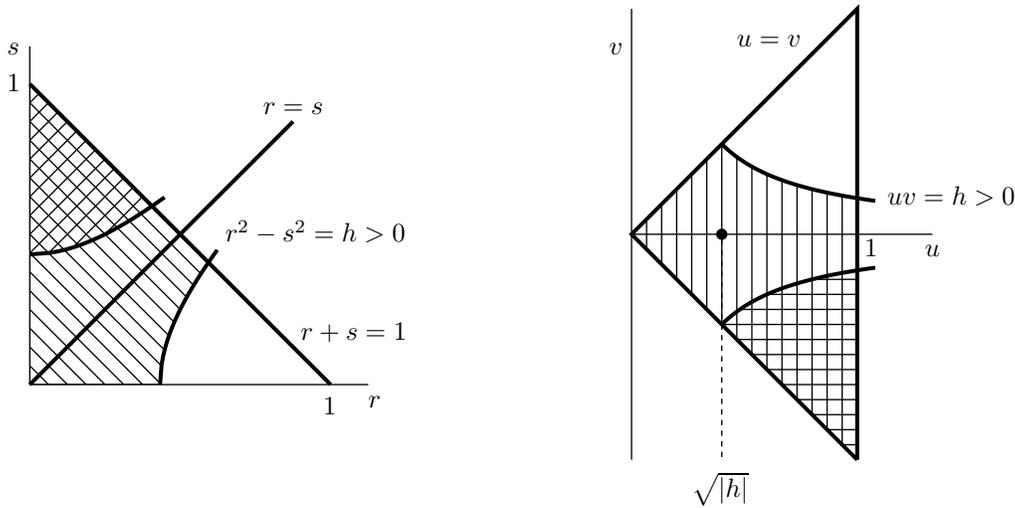}}
\put( 65, 2){$r$}
\put( 17,49){$s$}
\put( 59, 1){$1$}
\put( 17,44){$1$}
\put( 51,41){$r=s$}
\put( 56,11){$r+s=1$}
\put( 46,24){$r^2-s^2=h>0$}
\put(139,22){$u$}
\put( 97,49){$v$}
\put(131,22){$1$}
\put(114,50){$u=v$}
\put(134,29){$uv=h>0$}
\put(108,-10){$\sqrt{|h|}$}
\end{picture}
\caption{\textit{Regions of integration for coordinates $(r,s)$ and $(u,v)$.
Integral $I_{p,q}(h)$ for $h < 0$: doubly hatched region; $I_{p,q}(h)$ for
$h \geq 0$: doubly and singly hatched region.}\label{fig:intgebied}}
\end{figure}

We first recall the definition of $B_0$. Let $\calO$ be a
neighbourhood as in proposition \ref{pro:fnfs} such that $f$ can be
put into standard form c). Changing to cylinder coordinates $(r,
\varphi, s, \psi)$ and after an appropriate scaling we may assume that
$f(r, \varphi, s, \psi) = r^2 - s^2$ and $B = \{(r, \varphi, s, \psi)
\in \fR \times S^{p-1} \times \fR \times S^{q-1} \;|\; r+s \leq 1\}$ is
a subset of $\calO$. Then we define $B_0 = \{(r, \varphi, s, \psi) \in
B \;|\; 0 \leq r^2 - s^2 \leq h\}$ where $h \geq 0$. Now we wish to compute
\begin{displaymath}
\int_{B_0} 1 \,d\xi\,d\eta = \int_{B_0} r^{p-1} s^{q-1} g_p(\phi) g_q(\psi) dr\,ds\;d\phi\,d\psi = c_{p,q} I_{p,q}(h).
\end{displaymath}
The Jacobian of the transformation to cylinder coordinates is given by
$r^{p-1} s^{q-1} g_p(\phi) g_q(\psi)$. The integrand does not depend
on the angles so we split off the angular part and denote the
integrals by $a_p$ and $a_q$ where $a_m$ is the 'area' of the $m-1$
sphere. To facilitate the computations we make one further
transformation: $u=r+s$, $v=r-s$. In the constant $c_{p,q}$ we absorb
the constants $a_p$, $a_q$ and the Jacobian of the change of
coordinates $(r,s)$ to $(u,v)$. Furthermore we distinguish $h < 0$ and $h \geq 0$ and to simpify notation we write $k(u,v) = (u+v)^{p-1}(u-v)^{q-1}$. See figure \ref{fig:intgebied} for the regions of integration.
\begin{displaymath}
I_{p,q}(h) = \left\{
\begin{array}{lclr}
I^-_{p,q}(h) &=& \int_{\sqrt{-h}}^1\int_{-u}^{h/u} k(u,v) dv\,du, & h < 0\\
I^+_{p,q}(h) &=& \int_0^1 \int_{-u}^0 k(u,v) dv\,du +
                 \int_0^{\sqrt{h}} \int_0^u k(u,v) dv\,du +
                 \int_{\sqrt{h}}^1 \int_0^{h/u} k(u,v) dv\,du, & h \geq 0
\end{array}\right.
\end{displaymath}

After some computations it turns out that $I_{p,q}(h)$ consists of
several parts: for both $h < 0$ and $h \geq 0$:
\begin{displaymath}
I^{\pm}_{p,q}(h) = P^{\pm}_{p,q}(h) \;+\; \alpha^{\pm}_{p,q}(\pm h)^{(p+q)/2}
 \;+\; \beta^{\pm}_{p,q}(\pm h)^{(p+q)/2} \log\sqrt{\pm h} \;+\; \gamma^{\pm}_{p,q}.
\end{displaymath}
There is a constant part $\gamma^{\pm}_{p,q}$ because we consider the
volume between levels $-1$ and $h$. This is done for computational
reasons. Then continuity demands that $\gamma^+_{p,q} =
\gamma^-_{p,q}$ which turns out to be true. Common to all cases is a
polynomial part $P^{\pm}_{p,q}(h)$. The part with coefficient
$\alpha^{\pm}_{p,q}$ may contain a square root depending on $p$ and
$q$ and if the coefficient $\beta^{\pm}_{p,q}$ is nonzero there is a
logarithmic part. The definitions of $P$, $\beta$, $\gamma$ and
$\alpha$ are as follows
\begin{eqnarray*}
P^+_{p,q}(h) = P^-_{p,q}(h) &=& \sum_{k,m} \strut^{(1)} \frac{\binom{p-1}{k}\binom{q-1}{m} (-1)^m h^{k+m+1}}{(k+m+1)(p+q-2(k+m+1))},\\
\beta^+_{p,q} = \beta^-_{p,q} &=& \sum_{k,m} \strut^{(2)} \frac{\binom{p-1}{k}\binom{q-1}{m} (-1)^m}{(p+q)},\\
\gamma^+_{p,q} = \gamma^-_{p,q} &=& \sum_{k,m} \frac{\binom{p-1}{k}\binom{q-1}{m} (-1)^{p-1-k}}{(p+q)(p+q-(k+m+1))},\\
\alpha^+_{p,q} = -\alpha^-_{q,p} &=& \sum_{k,m} \strut^{(1)} 2\frac{\binom{p-1}{k}\binom{q-1}{m} (-1)^m}{(p+q)(p+q-2(k+m+1))} + \sum_{k,m} \strut^{(2)} 2\frac{\binom{p-1}{k}\binom{q-1}{m} (-1)^m}{(p+q)^2}\\
 &=& \frac{2}{(p+q)} \sigma_1(p,q) + \frac{2}{(p+q)^2} \sigma_2(p,q).
\end{eqnarray*}
The sum $\sum_{k,m}^{(1)}$ is taken over all $k$ and $m$ satisfying $0
\leq k \leq p-1$, $0 \leq m \leq q-1$ and $2(k+m+1) \neq (p+q)$,
whereas the sum $\sum_{k,m}^{(2)}$ is taken over the same range of $k$
and $m$ but now $2(k+m+1)=(p+q)$. The last line defines $\sigma_1$ and
$\sigma_2$. With these definitions $\beta^+_{p,q} = \frac{2}{(p+q)}
\sigma_2(p,q)$. The following properties of $\sigma_1$ and $\sigma_2$
are easily checked.
\begin{enumerate}
\item $\sigma_1(p,q) = (-1)^q \sigma_1(p,q)$, so $\sigma_1(p,q)=0$ for $q$ odd,
\item $\sigma_2(p,q) = (-1)^{q-1} \sigma_2(p,q)$, so $\sigma_2(p,q)=0$ for $q$ even,
\item $\sigma_2(p,q) \neq 0$ only if $p+q$ even, so using 2) $\sigma_2(p,q) \neq 0$ only if both $p$ and $q$ are odd.
\end{enumerate}

Now the next lemma is immediate.

\begin{lemma}\label{lem:nondif}
For each $p$ and $q$ the $\lceil\frac{n}{2}\rceil$-th derivative of
$I_{p,q}$ as a function of $h$ is discontinuous at $h=0$. The nature
of the discontinuity depends on $p$ and $q$. For both $p$ and $q$ even
it is a jump, for $p$ and $q$ odd the discontinuity is a log-like
singularity and for $p+q$ odd it is a root-like singularity.
\end{lemma}

\subsection{The volume of $B_0$ in proposition \ref{pro:relzad}}\label{app:borelzad}

The computation of the volume of $B_0$ in proposition \ref{pro:relzad} is similar to that in section \ref{app:bozad}. Here we only indicate the differences.

First note that $p$ and $q$ have a slightly different meaning because $p+q=n-1$. In this case the function $k$ in the expression for $I_{p,q}(h)$ is given by $k(u,v) = (h-uv)(u+v)^{p-1}(u-v)^{q-1}$. With this definition of $k$ the functions $I^{\pm}_{p,q}(h)$ are defined as before. Again after some computations we find
\begin{displaymath}
I^{\pm}_{p,q}(h) = P^{\pm}_{p,q}(h) + \alpha^{\pm}_{p,q}(\pm h)^{(p+q+2)/2} + \beta^{\pm}_{p,q}(\pm h)^{(p+q+2)/2} \log\sqrt{\pm h} + \gamma^{\pm}_{p,q} + \delta^{\pm}_{p,q}h.
\end{displaymath}
The expressions for $P$, $\alpha$, $\beta$ and $\gamma$ are more involved than in the previous section. Their structure, however, is similar therefore we skip the details. The relations are equal. There is one new term in the expression above which is defined as
\begin{eqnarray*}
\delta^+_{p,q} = \delta^-_{p,q} &=& \sum_{k,m} \frac{\binom{p-1}{k}\binom{q-1}{m} (-1)^{p-1-k}}{(p+q)(p+q-(k+m+1))}.
\end{eqnarray*}

Thus we come to the same conclusion as in section \ref{app:bozad}.

\begin{lemma}\label{lem:nondif2}
For each $p$ and $q$ the $\lceil\frac{n+1}{2}\rceil$-st derivative of
$I_{p,q}$ as a function of $h$ is discontinuous at $h=0$. The nature
of the discontinuity depends on $p$ and $q$. For both $p$ and $q$ even
it is a jump, for $p$ and $q$ odd the discontinuity is a log-like
singularity and for $p+q$ odd it is a root-like singularity.
\end{lemma}


\end{document}